\documentclass[oneside,article]{memoir}

\usepackage{amsmath}         
\usepackage{amsthm,thmtools,thm-restate} 
\usepackage{enumitem}        

\usepackage{rotating}

\usepackage[final]{listings}
\usepackage{bbm}
\usepackage{xcolor}

\definecolor{codegreen}{rgb}{0,0.6,0}
\definecolor{codegray}{rgb}{0.5,0.5,0.5}
\definecolor{codepurple}{rgb}{0.58,0,0.82}
\definecolor{backcolour}{rgb}{0.95,0.95,0.92}

\lstdefinestyle{mystyle}{
    backgroundcolor=\color{backcolour},   
    commentstyle=\color{codegreen},
    keywordstyle=\color{magenta},
    numberstyle=\tiny\color{codegray},
    stringstyle=\color{codepurple},
    basicstyle=\ttfamily\footnotesize,
    breakatwhitespace=false,         
    breaklines=true,                 
    captionpos=b,                    
    keepspaces=true,                 
    numbers=left,                    
    numbersep=5pt,                  
    showspaces=false,                
    showstringspaces=false,
    showtabs=false,                  
    tabsize=2
}

\usepackage{algorithm}
\usepackage{algpseudocode}
\usepackage{algorithmicx}
\usepackage{threeparttable}
\usepackage{float}

\usepackage{tikz}
\usetikzlibrary{matrix,arrows,positioning}
\usepackage{amsfonts}
\usepackage{parskip}
\usepackage{mathtools}
\usepackage{wrapfig}

\usepackage{tikz-cd}

\usepackage{caption, subcaption} 

\usepackage[hidelinks]{hyperref}                   
\usepackage[capitalize,nosort]{cleveref}           

\usepackage{amssymb}
\usepackage[mathscr]{euscript}
\usepackage{relsize}
\usepackage{stmaryrd}
\usepackage{stackengine}

\usepackage{indentfirst} 
\usepackage{multicol}

\usepackage{graphicx}
\usepackage{ragged2e} 
\usepackage[inline,nomargin]{fixme} 
\fxsetup{targetlayout=color}

\isopage[12]

\setlrmargins{*}{*}{1}
\checkandfixthelayout

\counterwithout{section}{chapter}
\usepackage{appendix}

\makeatletter
\providecommand{\abx@aux@refcontext}[1]{}
\providecommand{\abx@aux@cite}[2]{}
\providecommand{\abx@aux@segm}[3]{}
\providecommand{\abx@aux@page}[2]{}
\providecommand{\abx@aux@fnpage}[2]{}
\providecommand{\abx@aux@backref}[5]{}
\providecommand{\abx@aux@defaultrefcontext}[3]{}
\providecommand{\abx@aux@read@bbl@mdfivesum}[1]{}
\providecommand{\abx@aux@read@bblrerun}{}
\providecommand{\abx@aux@sortscheme}[1]{}
\providecommand{\abx@aux@refsection}[1]{}
\providecommand{\abx@aux@number}[2]{}
\makeatother

\makeatletter
\let\dgm\@undefined
\makeatother

  \declaretheorem[style=definition,within=section]{definition}
  \declaretheorem[style=definition,numberlike=definition]{example}

  \declaretheorem[style=plain,numbered=no,name=Theorem]{theorem*}

  \Crefname{corollary}{Corollary}{Corollaries}
  \Crefname{definition}{Definition}{Definitions}
  \Crefname{lemma}{Lemma}{Lemmas}
  \Crefname{proposition}{Proposition}{Propositions}
  \Crefname{remark}{Remark}{Remarks}
  \Crefname{theorem}{Theorem}{Theorems}
  \Crefname{notation}{Notation}{Notations}
  \Crefname{conjecture}{Conjecture}{Conjectures}

  \newlist{axioms}{enumerate}{1}
  \Crefname{axiomsi}{}{}

  \newenvironment{tikzeq*}
  {
    \begingroup
    \begin{equation*}
    \begin{tikzpicture}[baseline=(current bounding box.center)]
  }
  {
    \end{tikzpicture}
    \end{equation*}
    \endgroup
    \ignorespacesafterend
  }

  \tikzset
  {
    diagram/.style=
    {
      matrix of math nodes,
      column sep={4.3em,between origins},
      row sep={4em,between origins},
      text height=1.5ex,
      text depth=.25ex
    },
    over/.style={preaction={draw=white,-,line width=6pt}},
    every to/.style={font=\footnotesize},
    inj/.style={right hook->},
    surj/.style={-{Latex[open]}},
    cof/.style={>->},
    fib/.style={->>},
  }

  \DeclareFontFamily{U}{mathx}{\hyphenchar\font45}

  \DeclareFontShape{U}{mathx}{m}{n}{
    <5> <6> <7> <8> <9> <10>
    <10.95> <12> <14.4> <17.28> <20.74> <24.88>
    mathx10}{}

  \DeclareSymbolFont{mathx}{U}{mathx}{m}{n}

  \DeclareFontFamily{U}{mathb}{\hyphenchar\font45}

  \DeclareFontShape{U}{mathb}{m}{n}{
    <5> <6> <7> <8> <9> <10>
    <10.95> <12> <14.4> <17.28> <20.74> <24.88>
    mathb10}{}

  \DeclareSymbolFont{mathb}{U}{mathb}{m}{n}

  \DeclareMathAccent{\widebar}{0}{mathx}{"73}

  \DeclareMathSymbol{\Rsh}{\mathrel}{mathb}{"E9}

  \DeclareFontFamily{U}{MnSymbolA}{}

  \DeclareFontShape{U}{MnSymbolA}{m}{n}{
    <-6> MnSymbolA5
    <6-7> MnSymbolA6
    <7-8> MnSymbolA7
    <8-9> MnSymbolA8
    <9-10> MnSymbolA9
    <10-12> MnSymbolA10
    <12-> MnSymbolA12}{}

  \DeclareSymbolFont{MnSyA}{U}{MnSymbolA}{m}{n}

  \DeclareMathSymbol{\twoheaddownarrow}{\mathrel}{MnSyA}{27}

  \newcommand{\MSC}[1]{%
    \let\thempfn\relax
    \footnotetext[0]{2020 Mathematics Subject Classification: #1.}
  }

  \newcommand{\VR}{\textsf{VR}}
  \newcommand{\dgm}{\textsf{dgm}}
  \newcommand{\DCH}{\textsf{DH}}

\tikzstyle{vertex}=[circle, draw, minimum size=7pt, inner sep=0pt]

\DeclareFontFamily{U}{dmjhira}{}
\DeclareFontShape{U}{dmjhira}{m}{n}{ <-> dmjhira }{}

\author{Krzysztof Kapulkin \and Nathan Kershaw}

\title{Topological data analysis using persistent discrete homology}

\date{\today}

\begin{document}

\maketitle

\begin{abstract}
  We propose persistent discrete homology as a tool for topological data analysis and discuss its advantages over the existing methods.
  In particular, we provide empirical evidence that persistent discrete homology is more noise-resistant than persistent homology of the Vietoris--Rips complex for data coming from non-metric settings.

\end{abstract}



\section*{Introduction}

Since its introduction \cite{zomorodian-carlsson,carlsson:topology-and-data}, topological data analysis has seen a wealth of applications, spanning areas as disparate as time series analysis \cite{royer2021atol,umeda2017time}, neuroscience \cite{gardner-et-al,schneider-lee-mathis}, biology \cite{benjamin-et-al}, medicine \cite{dindin2020topological}, and physics \cite{sale-lucini-giansiracusa}; and this list is far from exhaustive.
Its essential premise is that many data sets encountered in practice carry interesting topological features, and understanding these features aids the analysis of the given data set.
As topological features of a space are usually detected using homology groups, persistent homology has become a key tool of topological data analysis.
Persistent homology takes a data set to a filtration of simplicial/cubical/chain complexes and applies homology to each element of the filtration.
By tracing the generators through different scales, we can then identify features that \emph{persist} (hence `persistent' in `persistent homology').

Several filtrations have been considered in the field, e.g., Delaunay or \v{C}ech filtrations \cite{edelsbrunner-mucke-alpha-shapes, bauer-edelsbrunner-cech-delaunay}, with the Vietoris--Rips filtration \cite{edelsbrunner-letscher-zomorodian,zomorodian-carlsson,edelsbrunner-harer:computational-topology} receiving the most attention.
Despite its success, a common drawback of Vietoris--Rips persistence is its lack of resistance to noise and outliers, and many techniques have been introduced to address this problem, amongst which multiparameter persistent homology has been studied most extensively \cite{carlsson-zomorodian-multiperam, Botnan-Lesnick-multiperam-intro, blumberg-lesnick-2-peram-stability}.
While these approaches have seen spectacular success in the metric setting (i.e., when the data is a metric space), in a non-metric setting a new problem arises.
Data coming from semi-metric settings is important to address: the semi-distances can represent any notion of similarity not satisfying the triangle inequality, such as correlations between two time series or blocks of text.
In such settings, the triangle inequality is precisely the constraint that ordinarily forces `short cycles' to be filled soon after they form.
When the inequality fails, this control is lost: a short cycle can remain unfilled across a wide range of scales.
Since these short cycles are easily created by noise, generators caused by noise persist longer and become much harder to distinguish from the actual signal.

The main purpose of this paper is to introduce persistent discrete homology as a new technique to address the non-metric setting.
Discrete homology was developed as a homology theory of graphs \cite{barcelo-capraro-white} and studied extensively from the combinatorial perspective \cite{barcelo-greene-jarrah-welker:comparison,barcelo-greene-jarrah-welker:vanishing,carranza-kapulkin:cubical-graphs} thanks to its broad applications in areas such as matroid theory and hyperplane arrangements \cite{barcelo-laubenbacher}.
Our key insight is that the same tool can be repurposed for data analysis, and that it addresses the non-metric problem at its source.
The difference between discrete homology and the simplicial homology of the clique complex, as used in Vietoris--Rips persistence, lies in which cycles it ``fills.''
In degree 1, both compute the homology of a 2-dimensional cell complex built out of a graph, and for both the 1-skeleton is the graph itself; but where Vietoris--Rips adds a 2-cell for each 3-cycle, discrete homology adds a 2-cell for every 3- \emph{and} 4-cycle.
This distinction is exactly what matters in the non-metric setting.
The 4-cycles that go unfilled by Vietoris--Rips --- and that, as noted above, can persist across a wide range of scales once the triangle inequality fails --- are precisely the ones most easily created by noise, and they generally carry little topological signal.
Discrete homology fills them by construction, so these spurious generators never appear in persistent discrete homology at all, while the theory retains the genuine signal.
This resistance is not merely empirical: stability of persistent discrete homology, established by Bubenik and Mili\'{c}evi\'{c} in the setting of closure spaces \cite{bubenik-milicevic}, applies directly in our context and provides a theoretical footing for the technique.

To quantify resistance to noise, we introduce a notion of similarity, called the overlap, between two persistence diagrams, i.e., outputs of persistent homology.
By analyzing the overlap between diagrams of data with and without noise, we can ascertain how well a given tool retains the topological signal.
We find that while in the metric setting Vietoris--Rips persistence behaves identically to discrete homology, in the non-metric setting, discrete homology retains significantly more signal.
This naturally raises the question of whether ``filling in'' higher cycles may result in a homology theory that retains more topological signal under noise.
The answer to this question, perhaps surprisingly, is negative.
As we start filling in 5- or 6-cycles, we notice that while the method sees less noise, it also fails to see meaningful topological signal.

In addition to testing our technique on curated data sets with known topological properties, we also test it on weather data, one obtained through a simulation and one from the U.S.~National Oceanic and Atmospheric Administration \cite{noaa}.
We are able to show that persistent discrete homology is better at detecting features within these data and distinguishing them from noise.

\textbf{Related work.}
The connection between discrete homotopy theory \cite{babson-barcelo-longueville-laubenbacher,barcelo-laubenbacher,carranza-kapulkin:cubical-graphs}, the general area studying invariants like discrete cubical homology, and data analysis was briefly discussed in \cite{memoli-zhou}.
Bubenik and Mili\'{c}evi\'{c} proved stability of persistent discrete homology in the context of closure spaces in \cite{bubenik-milicevic}, which directly applies in our setting.

\textbf{Organization of the paper.}
The first two sections of the paper review the necessary background on persistent homology (\cref{sec:pers_homology}) and discrete homology respectively (\cref{sec:discrete_homology}).
In \cref{sec:DCH_vs_TDA}, we introduce the notion of overlap of persistence diagrams and use it to compare the two homology theories.
Subsequently in \cref{sec:data}, we show that persistent discrete homology is more effective at detecting features in weather-related data than persistent homology of the Vietoris--Rips filtration.
In \cref{sec:summary}, we summarize the contributions of the present paper and discuss ideas for future work.
\section{Persistent homology and the Vietoris--Rips filtration} \label{sec:pers_homology}

In this section, we briefly review the necessary background on persistent homology.
Details can be found in a variety of excellent sources, including: \cite{zomorodian-carlsson,carlsson:topology-and-data, oudot-tda}.

In the most general case, a dataset can be viewed as a finite semi-metric space $(X,d)$, where $d$ might fail the triangle inequality, and $d(x,y)=0$ need not imply $x=y$.

\begin{definition} \label{def:rips_cpx}
    Let $(X,d)$ be a finite semi-metric space, and $r \in \mathbb{R}_{\geq 0}$. 
    The \emph{Vietoris--Rips complex at scale $r$} is the simplicial complex $\textsf{VR}_r(X,d)$, which has the vertex set $\textsf{VR}_r(X,d)_0 = X$, and an $n$-simplex $\{x_0, \dots, x_n\} \in \textsf{VR}_r(X,d)_n$ if $d(x_i,x_j) \leq r$ for all $0 \leq i,j \leq n$.
\end{definition}

By letting the value of $r$ vary, we obtain a filtration of simplicial complexes, called the Vietoris--Rips filtration, with natural inclusions $\VR_r(X,d) \hookrightarrow \VR_s(X,d)$.
Taking homology $\mathcal{H}_n(-,\mathbb{F})$ with values in a field $\mathbb{F}$ of this filtration yields a sequence of $\mathbb{F}$-vector spaces $\mathcal{H}_n(\VR(X,d),\mathbb{F})$, where $\mathcal{H}_n(\VR(X,d),\mathbb{F})_r = \mathcal{H}_n(\VR_r(X,d),\mathbb{F})$.

Using a decomposition theorem \cite{gabriel:decomposition}, we can convert this sequence of vector spaces into a \emph{persistence diagram}, denoted $\dgm_n(\VR(X,d))$, which is a multiset of pairs $(b,d)$, with each pair $(b,d)$ corresponding to a homology generator born at scale $b$ and dying at scale $d$.

To visualize a persistence diagram, we will draw its \emph{barcode}, which plots horizontal bars from $b$ to $d$ for every pair $(b,d)$ in the persistence diagram. 
An example of this can be seen in:
\begin{figure}[H]
    \centering
    \begin{subfigure}[b]{0.38\textwidth}
        \centering
        \includegraphics[width=\textwidth]{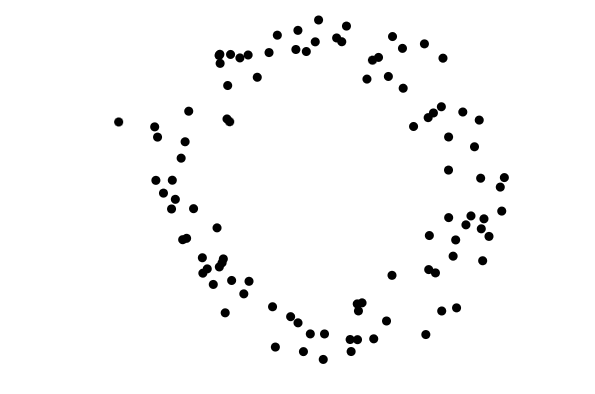}
    \end{subfigure}
    \hfill
    \begin{subfigure}[b]{0.58\textwidth}
        \centering
        \includegraphics[width=\textwidth]{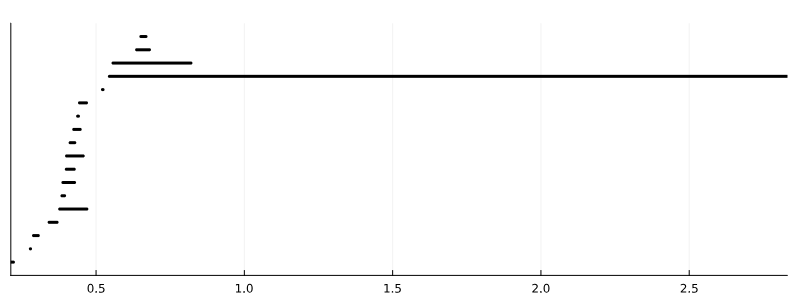}
    \end{subfigure}
    \caption{Sample data (left) and dimension one barcode (right)}
\end{figure}
Here, the long bar in the barcode represents the generator coming from the overall circle-shape of the data. 
The shorter bars come from generators that appear and disappear relatively quickly, which are due to the noise in the data. 
Using this pipeline, we are able to convert our data into these persistence diagrams (or barcodes).
The longer generators generally inform us on the larger topological features in our data, and the shorter generators generally arise from noise. 

\section{Discrete homology} \label{sec:discrete_homology}
Discrete homology was introduced in \cite{barcelo-capraro-white} in the context of finite metric spaces, and studied in the context of graph theory in a series of papers by Barcelo, Greene, Jarrah, and Welker \cite{barcelo-greene-jarrah-welker:comparison,barcelo-greene-jarrah-welker:vanishing,barcelo-greene-jarrah-welker,carranza-kapulkin:cubical-graphs}.
Much work has also been done on computability \cite{kapulkin-kershaw, ender-kapulkin-first, ender-kapulkin-higher}, rendering it a viable tool for topological data analysis.

For us, a \emph{graph} $G$ consists of a set with a symmetric and reflexive relation.
We write $G_V$ for the underlying set, called the \emph{vertex set}, and $G_E \subseteq G_V \times G_V$ for the relation, which we call the \emph{edge set}.
We write $v \sim w$ to denote $(v,w) \in G_E$.
A \emph{graph map} $f \colon G \to H$ from a graph $G$ to $H$ is a set map $f \colon G_V \to H_V$ that preserves the relation.
Put differently, our graphs are undirected, simple (no multiple edges), and reflexive (each vertex has a unique loop).
Graph maps take vertices to vertices preserving the adjacency (edge) relation, but since the graphs are reflexive, the graph maps are in particular allowed to ``contract'' edges.

Several examples of graphs are important in our setting.
The \emph{line graph} $I_n$ has $n+1$ vertices, labelled $0, \ldots, n$ with an edge between vertices $i$, $i+1$ for all $0 \leq i \leq n-1$; the $n$-\emph{cycle graph} $C_n$ is then the quotient of $I_n$ by identifying $0$ and $n$.
Given graphs $G$ and $H$, their \emph{box product} $G \square H$ has the vertex set $G_V \times H_V$, and an edge between vertices $(v,w)$ and $(v',w')$ if either $v=v'$ and $w \sim w'$, or $v \sim v'$ and $w = w'$. 
The \emph{discrete $n$-cube} (or the \emph{hypercube graph}), denoted $I_1^{\square n}$ is the iterated box product of $I_1$ with itself $n$ times.
We label the vertices $(x_1,\dots,x_n)$ for $x_i \in \{0,1\}$.
\begin{figure}[H]
  \centering
  \begin{subfigure}[b]{0.19\textwidth}
    \centering
    \begin{tikzpicture}[scale=.9]
      \foreach \x/\lab in {0/0,1/1,2/2,3/3}
        \node[circle,draw,inner sep=2pt,label=above:\lab] (\x) at (\x,0) {};
      \draw (0)--(1)--(2)--(3);
    \end{tikzpicture}
    \caption{$I_{3}$}\label{fig:sub1}
  \end{subfigure}\hfill
  \begin{subfigure}[b]{0.19\textwidth}
    \centering
    \begin{tikzpicture}[scale=.9]
      \node[circle,draw,inner sep=2pt,label=above:0] (0) at ( 0,   0.78) {};
      \node[circle,draw,inner sep=2pt,label=right:1] (1) at ( 0.95,1.31) {};
      \node[circle,draw,inner sep=2pt,label=above:2] (2) at ( 0.59,2.19) {};
      \node[circle,draw,inner sep=2pt,label=above:3] (3) at (-0.59,2.19) {};
      \node[circle,draw,inner sep=2pt,label=left:4]  (4) at (-0.95,1.31) {};
      \draw (0)--(1)--(2)--(3)--(4)--(0);
    \end{tikzpicture}
    \caption{$C_{5}$}\label{fig:sub2}
  \end{subfigure}\hfill
  \begin{subfigure}[b]{0.19\textwidth}
    \centering
    \begin{tikzpicture}
      \node[circle,draw,inner sep=2pt] (0) at (0,0) {};
      \node[circle,draw,inner sep=2pt] (1) at (1,0) {};
      \draw (0)--(1);
    \end{tikzpicture}
    \caption{$I_1$}\label{fig:I3}
  \end{subfigure}\hfill
  \begin{subfigure}[b]{0.19\textwidth}
    \centering
    \begin{tikzpicture}
      \node[circle,draw,inner sep=2pt] (0) at (0,0) {};
      \node[circle,draw,inner sep=2pt] (1) at (0,1) {};
      \node[circle,draw,inner sep=2pt] (2) at (1,0) {};
      \node[circle,draw,inner sep=2pt] (3) at (1,1) {};
      \draw (0)--(1);
      \draw (1)--(3);
      \draw (2)--(3);
      \draw (2)--(0);
    \end{tikzpicture}
    \caption{$I_1^{\square 2}$}\label{fig:C5_1}
  \end{subfigure}\hfill
  \begin{subfigure}[b]{0.19\textwidth}
    \centering
    \begin{tikzpicture}[scale=.8]
      \node[circle,draw,inner sep=2pt] (0) at (0,0) {};
      \node[circle,draw,inner sep=2pt] (1) at (0,1) {};
      \node[circle,draw,inner sep=2pt] (2) at (1,0) {};
      \node[circle,draw,inner sep=2pt] (3) at (1,1) {};
      \node[circle,draw,inner sep=2pt] (4) at (0.5,0.5) {};
      \node[circle,draw,inner sep=2pt] (5) at (0.5,1.5) {};
      \node[circle,draw,inner sep=2pt] (6) at (1.5,0.5) {};
      \node[circle,draw,inner sep=2pt] (7) at (1.5,1.5) {};
      \draw (0)--(1); \draw (1)--(3); \draw (2)--(3); \draw (2)--(0);
      \draw (4)--(5); \draw (5)--(7); \draw (6)--(7); \draw (6)--(4);
      \draw (0)--(4); \draw (1)--(5); \draw (2)--(6); \draw (3)--(7);
    \end{tikzpicture}
    \caption{$I_1^{\square 3}$}\label{fig:C5_2}
  \end{subfigure}
  \caption{Examples of graphs.}
  \label{fig:box powers}
\end{figure}

\begin{definition} \label{def:sing_n-cube}  \leavevmode
  \begin{itemize}
      \item  A \emph{(singular) $n$-cube} in a graph $G$ is a graph map $A \colon I_1^{\square n} \to G$. 
      \item For a singular $n$-cube $A$ in a graph $G$ and $1 \leq i \leq n$, we define the \emph{positive} and \emph{negative face maps} $\delta_i^- A$ and $\delta_i^+ A$ to be the singular $(n-1)$-cubes given by:
    \[\delta_i^+ A (x_1, \dots, x_{n-1}) := A(x_1, \dots, x_{n-1}, 1, x_i, \dots, x_n)\]
    \[\delta_i^- A (x_1, \dots, x_{n-1}) := A(x_1, \dots, x_{n-1}, 0, x_i, \dots, x_n)\]
    \item     If $\delta_i^- A = \delta_i^+ A$ for some $i$, we say that $A$ is \emph{degenerate}.
    Otherwise, we say that $A$ is \emph{non-degenerate}.
    \item For $n \geq 0$, define $\mathcal{C}_n(G)$ to be the free $R$-module generated by the set of non-degenerate singular $n$-cubes in $G$.
    These assemble into a chain complex $\mathcal{C}(G) = (\mathcal{C}_\bullet, \partial_\bullet)$ by defining $\partial_n \colon \mathcal{L}_n(G) \to \mathcal{L}_{n-1}(G)$ on generators by 
    \[ \partial_n(A) = \sum_{i=1}^n(-1)^i(\delta_i^-A -\delta_i^+ A)\text{.}\]
    and extending linearly.
    \item  For $n \geq 0$, the \emph{$n$-th discrete homology group} of a graph $G$ with coefficients in $R$ is $\mathcal{H}_n(G; R)$ is the homology of the chain complex $\mathcal{C}(G)$.
  \end{itemize}   
\end{definition}

\begin{example} \label{ex:homology_groups}
    The first homology groups of the $n$-cycle $C_n$ are given by: 
    $$\mathcal{H}_1(C_n) = \begin{cases}
        0 & n\leq 4 \\
        R & n \geq 5
    \end{cases}$$
\end{example}
Moreover, all higher homology groups of $C_n$ vanish. 
Thus, discrete homology sees $C_n$ as the circle for $n \geq 5$, and as contractible for $n \leq 4$.

Having now defined discrete homology, we briefly describe how it can be used in place of Vietoris--Rips persistence for persistent homology.
Namely, we have the following filtration of graphs:

\begin{definition}
    Let $X=(X,d)$ be a finite metric space. For a value $r \in \mathbb{R}_{\geq 0}$, let $\DCH_r(X,d)$ be the graph with the vertex set $\DCH_r(X)_V = X$ and the edge relation $x \sim x'$ if $d(x,x') \leq r$.
    We use $\DCH(X)$ to denote the corresponding filtration of graphs.
\end{definition}

We can use the discrete homology functor to convert $\DCH(X)$ into a persistence module, whose decomposition gives us a persistence diagram. 
It is known that persistent discrete homology is stable with respect to the Gromov--Hausdorff and bottleneck distances \cite{bubenik-milicevic}.

\section{Persistent discrete homology vs the Vietoris-Rips filtration} \label{sec:DCH_vs_TDA}

Having now defined two different methods of analyzing a dataset, the purpose of this section is to compare the two methods on some sample data. 
The results of this section are for degree 1 homology only.
We believe analogous results should hold for higher dimensions as well, however at present, we lack the formal geometric interpretation and computational tools to efficiently analyze dimension 2 for discrete homology. 
All persistent homology computations were done using the GUDHI library in Python \cite{gudhi:urm, gudhi:FilteredComplexes}, and are all done over $\mathbb{Z}_2$.

The main focus of this section is noise resistance, i.e., how each method behaves when adding noise into the dataset. 
What we would like to test is how similar the persistence diagrams are prior to and after adding noise. 
In order to do this in a precise manner, we thus need a notion of similarity between barcodes. 

The bottleneck distance, while a natural choice, is insufficient for our purposes for two reasons.
First, since the cost of a matching is given by a supremum, it is blind to adding a lot of short bars, as long as one bar is changed by a slightly larger length.
Second, it is not scale-invariant, and persistent discrete homology tends to slightly scale the intervals down. 
This would create a bias towards persistent discrete homology, as the distances would be lower due to the smaller scales.
While $p$-Wasserstein distances could be used to solve the first problem, they do not address the second.
We remark here, however, that, if done carefully, a normalized version of a $p$-Wasserstein distance could be an alternative to the methods of this paper.
Instead, we develop our own notion of similarity that is well suited for this analysis. 

\begin{definition} \label{def:length}
    Let $P$ be a persistence diagram. The \textit{length} of $P$ is given by:
    $$\textsf{len}(P) = \sum_{(b,d) \in P} d-b$$
    That is, the sum of all the lengths of intervals in $P$. 
\end{definition}

\begin{definition}
    let $P$ and $Q$ be persistence diagrams.
    A \textit{matching} between $P$ and $Q$ is a subset $M \subseteq P \times Q$, such that the projections $\pi_P:M \to P$ and $\pi_Q:M \to Q$ are injective.
\end{definition}

\begin{definition} \label{def:intersection-matching}
    Let $P$ and $Q$ be two persistence diagrams, and $M$ a matching between them. 
    The \textit{intersection} of $P$ and $Q$ over $M$ is the persistence module:
    $$P \cap_MQ = \{(b,d) \cap(b',d') \ | \ ((b,d),(b',d')) \in M \}$$
    That is, the collection of intersections of the matched intervals.
\end{definition}

\begin{definition} \label{def:intersection-len}
    Let $P$ and $Q$ be persistence modules. The \textit{intersection length} of $P$ and $Q$ is:
    $$\textsf{len}(P\cap Q) = \sup \{ \textsf{len}(P \cap_M Q) \ | \ M \text{ is a matching between } P \text{ and } Q \}$$
\end{definition}

We can now define our notion of similarity:
\begin{definition} \label{def:overlap}
    Let $P$ and $Q$ be persistence modules. 
    The \textit{overlap} between $P$ and $Q$ is given by:

    $$O(P,Q) = \frac{\textsf{len}(P\cap Q)^2}{\textsf{len}(P)\textsf{len}(Q)}$$
\end{definition}

We have that $O(P,Q) \in [0,1]$, where $O(P,Q)=1$ means that $P=Q$, and $O(P,Q)=0$ means they are disjoint. 
Since we are normalizing by both the length of $P$ and the length of $Q$ this will also be invariant of scale. 
An informal way of thinking about this is the percent of the signal in $P$ that is present in $Q$.

We can now introduce our methodology for testing the noise resistance of persistent discrete homology and homology of the Vietoris-Rips filtration. 
We want to test noise resistance in two cases: one is the setting where the data comes from a metric space, and the other is the setting where the data comes from a semi-metric space (failing the triangle inequality). 

In both cases, we begin with a dataset $X$ of 25 evenly spaced points about a unit circle in $\mathbb{R}^2$. 
Given a standard deviation $\sigma$, we then add noise in one of two ways:
\begin{itemize}
    \item[1)] Metric preserving: for each point $x \in X$, we move $x$ by adding $(e_1,e_2)$ to $x$, where $(e_1,e_2) \sim$ $N(0,\sigma^2) \times N(0,\sigma^2)$
    \item[2)] Non-metric preserving: for each pair of points $x$ and $y$, we add $e$ to $d(x,y)$, where $e \sim N(0,\sigma^2)$.
    If this creates a negative distance, we set that distance to 0. 
\end{itemize}
Call this new dataset with noise added $X'$. 
For a given value of $\sigma$, the noise resistance can be measured via the expected value $\mathbb{E}[O(\dgm(X),\dgm(X'))]$, where here $\dgm(X)$ is meant to mean the persistence diagram using whichever method of persistent homology is being tested.

We can then vary $\sigma$ and run this 1000 times for each value of $\sigma$, method of persistence, and method of adding noise:
\begin{figure}[H]
    \centering
    \begin{subfigure}[b]{0.48\textwidth}
        \centering
        \includegraphics[width=\textwidth]{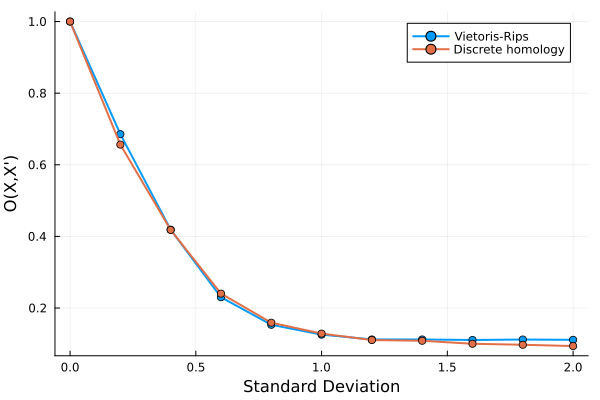}
    \end{subfigure}
    \hfill
    \begin{subfigure}[b]{0.48\textwidth}
        \centering
        \includegraphics[width=\textwidth]{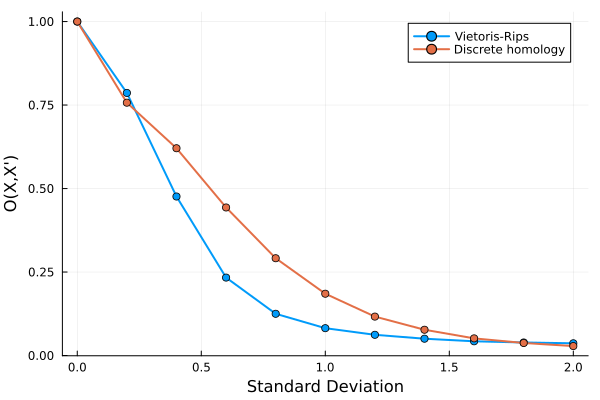}
    \end{subfigure}
    \caption{Noise resistance adding noise via method 1 (left) and method 2 (right)}
    \label{fig:noise_resistance}
\end{figure}
As one would expect, as the standard deviation increases, the noise resistance decreases. 
Of note is that using method 1, both discrete homology and the Vietoris-Rips filtration behave almost identically, meaning in terms of noise resistance, they are identical. 
Using method 2, however, discrete homology significantly outperforms Vietoris-Rips, specifically for data with a mild amount of noise added in. 
Thus, discrete homology appears to be more noise resistant for data equipped with a semimetric rather than a metric.

The reason why discrete homology is more noise resistant, specifically in the non-metric setting, is due to filling in both 3- and 4-cycles.
This is the main difference between discrete homology and homology of the Vietoris--Rips complex, which only fills in $3$-cycles.
At the same time, adding noise to data creates a large number of $4$-cycles.
Since the homology of the Vietoris--Rips complex will see the four cycles as holes and discrete homology will not, the barcode using the Vietoris--Rips filtration will always contain more bars than the barcode using discrete homology. 
In the metric setting, these bars corresponding to $4$-cycles will be relatively small, however, since due to the triangle inequality the diagonals of the square will appear in the filtration rather quickly. 
This is not the case in the non-metric setting, since the triangle inequality no longer holds. 
Since these four cycles are relatively easy to randomly create, this adds a lot of extra noise into the barcode in the non-metric setting. 

Since noise is created due to random $4$-cycles appearing, it follows as well that noise appears due to random $5$-cycles too, which is seen in both methods.
Thus, a natural idea would be to define a homology theory that fills in five cycles as well, or more generally, one filling in up to $n$ cycles for any $n$. 

One drawback is that those other homology theories do not as readily generalize to higher dimensions. 
It is also more computationally expensive the longer the cycles you are filling, due to both needing to identify those cycles in the graph, and also increasing the size of the corresponding matrix used in the persistence algorithm.
It may still be the case that these theories are more noise resistant; however. 

To test this, we can rerun the experiment using method two of adding noise, to test these different methods in the non-metric setting:
\begin{figure}[H]
    \centering
    \includegraphics[width=0.8\linewidth]{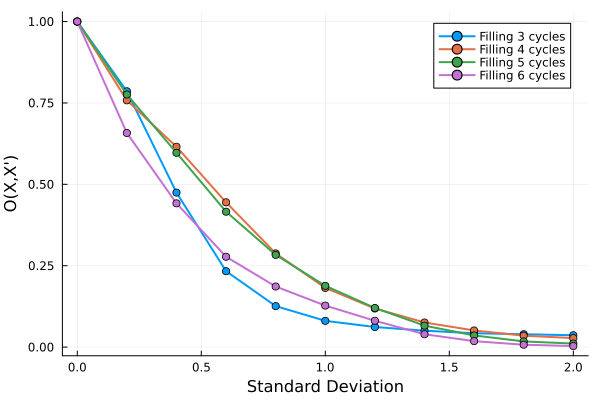}
    \caption{Noise resistance of different homology theories}
    \label{fig:all_cyc_lens}
\end{figure}
We see that the most noise-resistant methods are filling in up to 4-cycles (discrete homology) and filling in up to 5-cycles. 
Discrete homology is slightly better for some standard deviations, but the difference is not significant. 
Upon filling in 6-cycles as well, there is a steep drop-off in terms of noise resistance. 
Thus, for dimension one noise resistance, discrete homology and homology filling up to five cycles appear to be the best, however discrete homology is more practical due to being faster to compute and readily generalizing to higher dimensions. 

In light of the explanation for why discrete homology is more noise resistant, it may seem surprising that filling in up to 6-cycles is not more advantageous; however, this highlights the tradeoff one is making by filling in longer cycles. 
While it is true that you create less noise the longer cycles you fill, it also becomes easier to destroy the meaningful signal that was there originally. 
It turns out that for length 4, this tradeoff is advantageous, and for length 5, this tradeoff is relatively neutral, but for length 6, it starts to become a negative tradeoff. 
It is too hard to randomly create a 6-cycle with noise for the added reduction in noise to be worth the destroyed signal.

\section{Results on data} \label{sec:data}

In this section, we compare methods on more complex data.
Since we expect persistent discrete homology to perform better specifically for data equipped with a semi-metric, we focus on those types of data.
One such example, which is the primary focus of this section, is if the distance between points is given by some sort of a correlation. 
For example if our datapoints are weather stations, and the distance between two weather stations is dependent on the correlation coefficient between the readings between the two stations.
Two natural choices are to set $d(x,y) = 1-R^2$, or $d(x,y) = -\ln(R^2)$.
We focus on these types of distances.

\paragraph{Weather simulation.}
The first test we do is build and test a simulation meant to mimic weather-like data. 
The general setup is that there are several locations, each with a reading. 
These readings change over time with some amount of randomness. 
They not only change randomly, but also the readings at nearby locations influence the next reading at a given location. 
For example, the locations could be cities and the readings daily maximum temperatures. 
These will change randomly, such as due to clouds or wind, but the temperature at one location also influences temperatures at nearby locations (or at the very least, the same factors influence both). 

The setup is as follows: we fix an $n \times n$ grid, which we view as a graph $G = I_n \square I_n$. 
Each vertex represents a location. 
At time zero, assign each vertex $x$ a reading $y_x^0$ from some distribution, say, $y_x^0 \sim N(0,1)$. 
Given two vertices $w$ and $x$, the \textit{distance} between them, $d(w,x)$ is defined to be the length of the shortest path from $w$ to $x$.
Note that this distance is being used to generate readings; it is not the distance function on the dataset we use for persistence.
Given a vertex $x$, let $m_x = \textsf{max}\{d(x,w) \ | \ w \in G_V \}$. 
Given the collection of readings at a time $t$, and $x \in G_v$, we can define the readings at $t+1$ as follows:

\begin{enumerate}
    \item Given $x \in G_V$ let: 
    $$y_x = \frac{\Sigma_{z \in G_V} \ y_z^t(m_x - d(w,x) + 1)^2}{\Sigma_{z \in G_V} (m_x - d(z,x) + 1)^2}$$

    \item Then define ${y^t_x}' = y^t_x + e_x$, where $e_x\sim N(0,1)$

    \item Finally, set: $$y_x^{t+1} = \frac{\Sigma_{z \in G_V} \ {y^t_z}'(m_x - d(z,x) + 1)^2}{\Sigma_{z \in G_V} (m_x - d(z,x) + 1)^2}$$
\end{enumerate}

This allows us to get a time series of readings at each vertex, for any desired length.
The idea is that we start by using a weighted average inversely proportional to distance, and add some random noise. 
We then take a weighted average once more, so that the randomness introduced at each node also influences the nearby nodes. 

We then add in a disturbance to one of the vertices. 
To do this, we now view $G$ as a weighted graph with all edges weighted 1, and fix a weight $w \geq 1$. 
We then change the weights of the edges containing the chosen vertex to $w$. 
This influences the distance between nodes, as the distances are now computed using the weights. 
We can update the readings in the same manner.
The idea is that we will choose a random vertex in the interior of the grid, and we want to be able to detect which vertex was chosen.
In real-world scenarios, the vertex chosen represents an obstruction that influences the readings, for example, a lake or a mountain for weather data. 
We have three models for detecting the chosen vertex.

\emph{Expected vs actual model:}  
This most natural model compares the expected readings at time $t$ (given the $t-1$ readings and assuming all edges have a weight 1) and compares that to the actual readings.
We then sum the absolute value of the difference between the two over all readings, and interior vertex with the largest total difference is chosen.

\emph{Persistent discrete homology model:}
For each pair of vertices, we fit the readings of the two vertices against each other linearly. 
We then set $d(x,y)=1-R^2_{x,y}$. 
This gives us a semi-metric, allowing us to use persistent discrete homology.
We then look at the longest-lasting generator and find the vertices interior to it, using the representative of minimal cycle length. 
If only one vertex is interior, that one is chosen. 
If there are multiple interior vertices, we use the expected vs actual model to settle tiebreaks.

\emph{Vietoris--Rips model:} 
This is the same as the discrete homology model, but rather than taking persistent discrete homology, it uses the homology of the Vietoris--Rips filtration, and likewise defers to the expected vs actual model for tiebreaking.

We can set $G = I_8 \square I_8$, run with 50 readings on each node, and let $w$ vary, and check the accuracy of each model. 
We run with 5000 iterations for each value of $w$. 
The results can be seen in the figure below: 

\begin{figure}[H]
    \centering
    \includegraphics[width=0.8\linewidth]{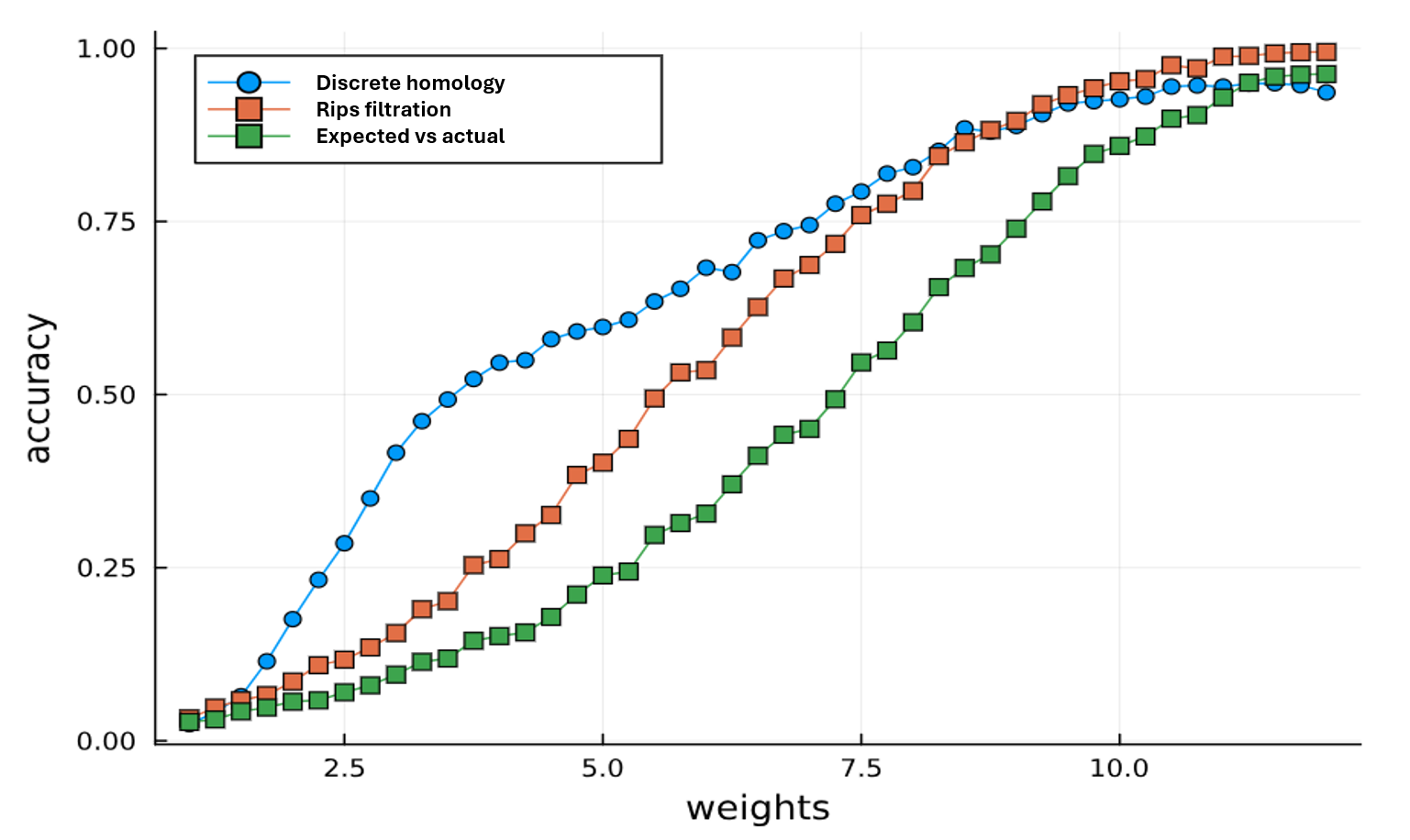}
    \caption{Accuracy of the three models}
    \label{fig:accuracy-of-modesl}
\end{figure}

Naturally, at very low values of $w$ (close to 1) the differences are negligible, and thus the chosen vertex is very difficult to detect; for high values of $w$, the differences are very large, and the chosen vertex is easy to detect with all models reaching accuracy close to $100\%$ for $w=12$. 
Importantly, however, the discrete homology model reaches a higher accuracy much faster than the other two. 
This suggests that the discrete homology model is better at handling noisier data, which is often the case with real-world data.

\paragraph{Weather data.} 
To show that the simulation discussed above is an accurate representation of reality, we use the global monthly summary data from the U.S.~National Oceanic and Atmospheric Administration (NOAA) \cite{noaa}. 
We run this for two regions: around Lake Ontario and around Lake Erie.
For Lake Ontario, we use a latitude range of (42.7,45) and a longitude range of (-80.5,-75). 
We then correlate the monthly precipitation readings between the stations. 
Since this often results in very high correlation coefficients, we choose to use the distance function of $d(x,y) = -\ln(R^2_{x,y})$. 
We can then use either persistent discrete homology or homology of the Vietoris--Rips filtration.
In both cases, the longest-lasting generator was the same. 
The weather stations, as well as the longest generator, can be seen below: 
\begin{figure}[H]
    \centering
    \begin{subfigure}[b]{0.48\textwidth}
        \centering
        \includegraphics[width=\textwidth]{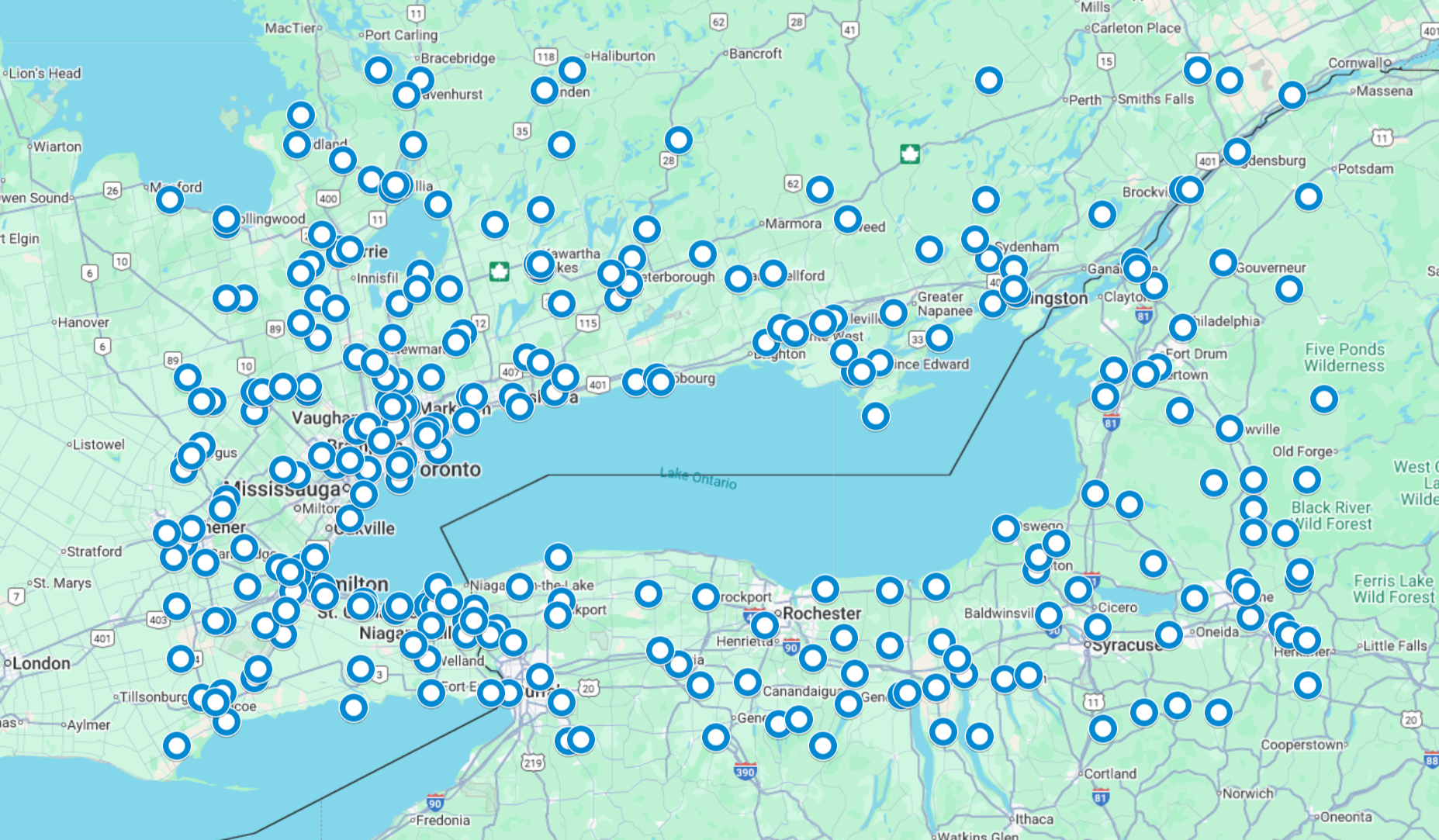}
    \end{subfigure}
    \hfill
    \begin{subfigure}[b]{0.48\textwidth}
        \centering
        \includegraphics[width=\textwidth]{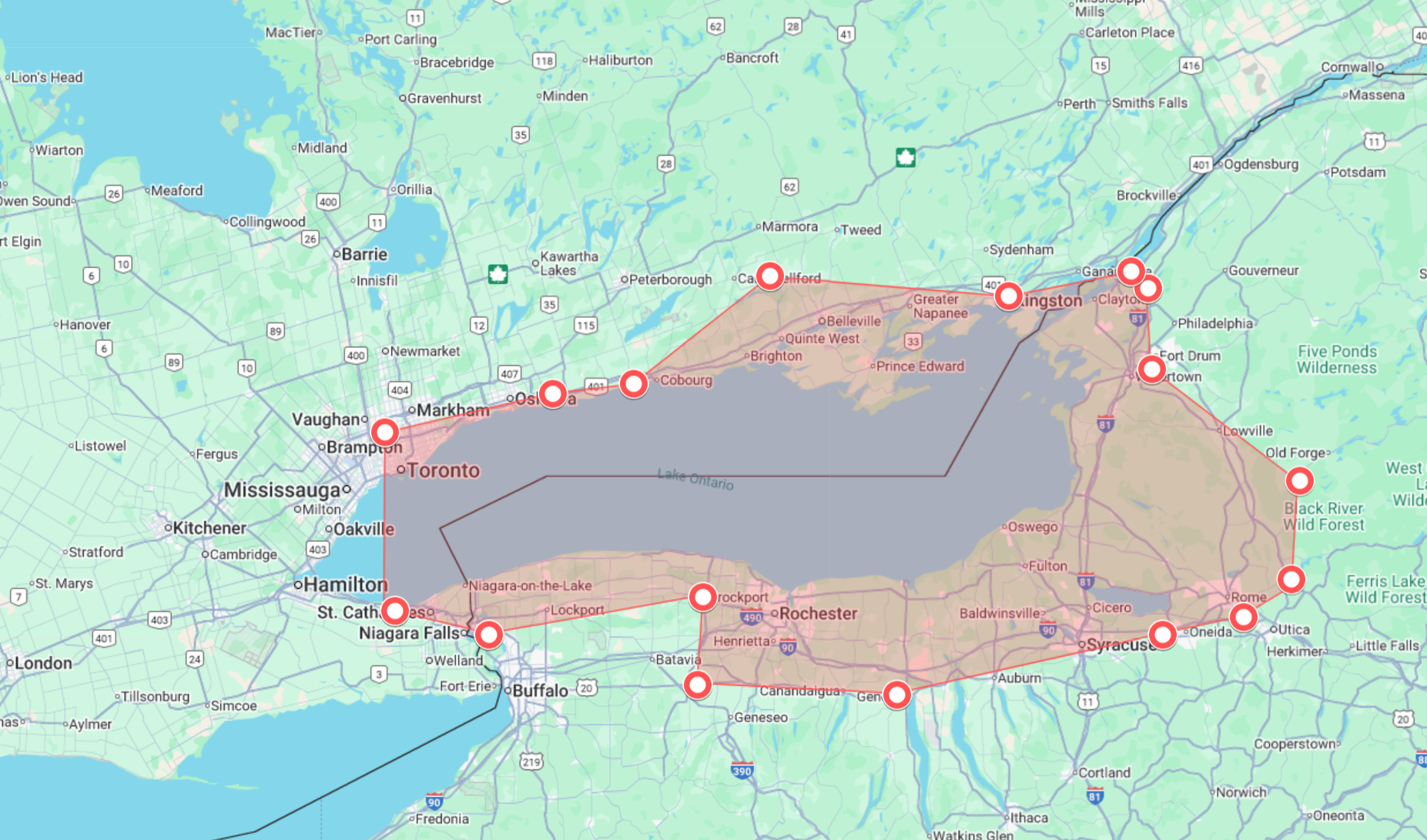}
    \end{subfigure}
    \caption{Lake Ontario weather stations (left) and longest generator (right)}
    \label{fig:station_and_cyc_ON}
\end{figure}
While both discrete homology and the Vietoris--Rips filtration detect the lake with the longest generator, their barcodes are significantly different: 
\begin{figure}[H]
    \centering
    \begin{subfigure}[b]{0.48\textwidth}
        \centering
        \includegraphics[width=\textwidth]{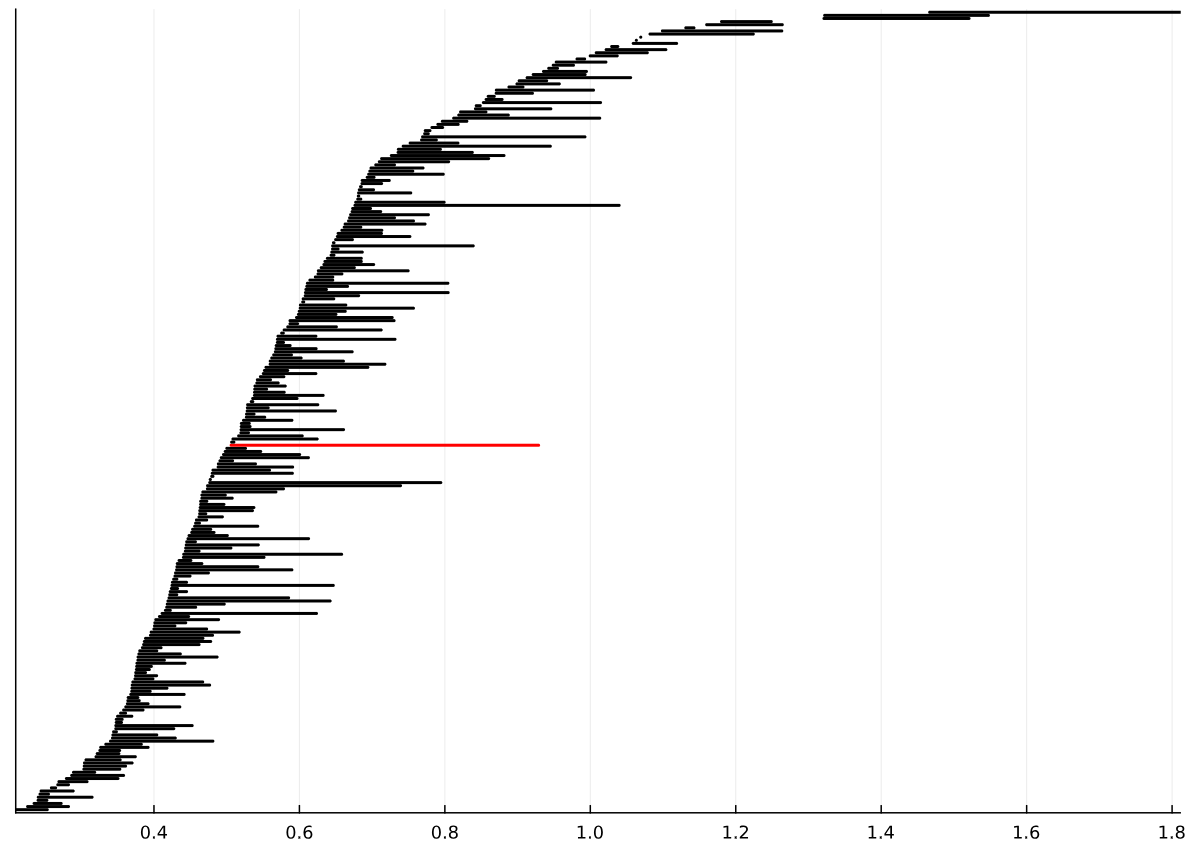}
    \end{subfigure}
    \hfill
    \begin{subfigure}[b]{0.48\textwidth}
        \centering
        \includegraphics[width=\textwidth]{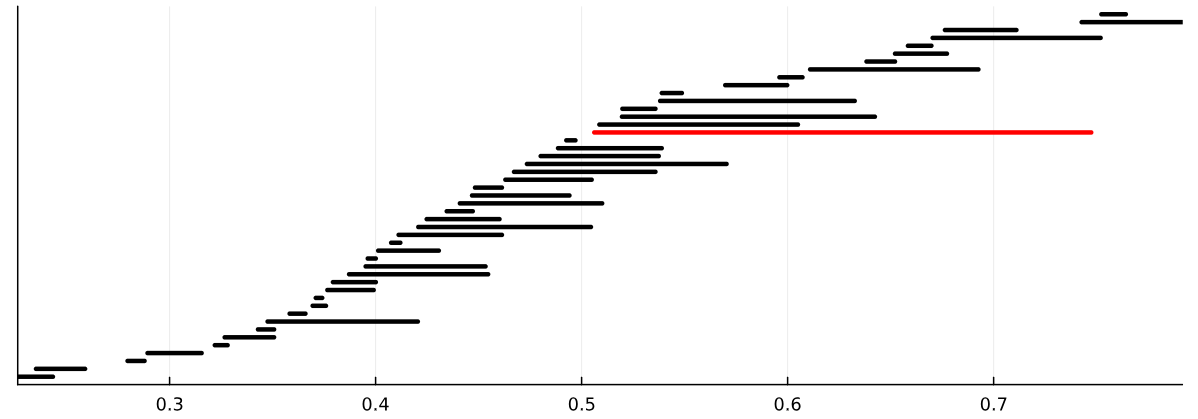}
    \end{subfigure}
    \caption{Barcodes using Vietoris-Rips (left) and discrete homology (right)}
    \label{fig:on_barcodes}
\end{figure}
The longest generator in each is highlighted in red. 
One would expect a large amount of noise regardless, due to the nature of detecting lakes using homology on correlations of precipitation data, but there is far more noise using Vietoris--Rips as opposed to discrete homology. 
For Vietoris--Rips, the longest bar is barely so, and almost indistinguishable from the other bars appearing due to noise. 
For discrete homology, while there is still a large amount of noise, the longest bar is almost double the length of the others.
It is much easier to differentiate from the noise. 

For Lake Erie, we use a latitude range of (41.1, 46.7) and a longitude range of (-89.9,-83.5). 
Again, the longest generator is the same for both, and can be seen below:

\begin{figure}[H]
    \centering
    \begin{subfigure}[b]{0.48\textwidth}
        \centering
        \includegraphics[width=\textwidth]{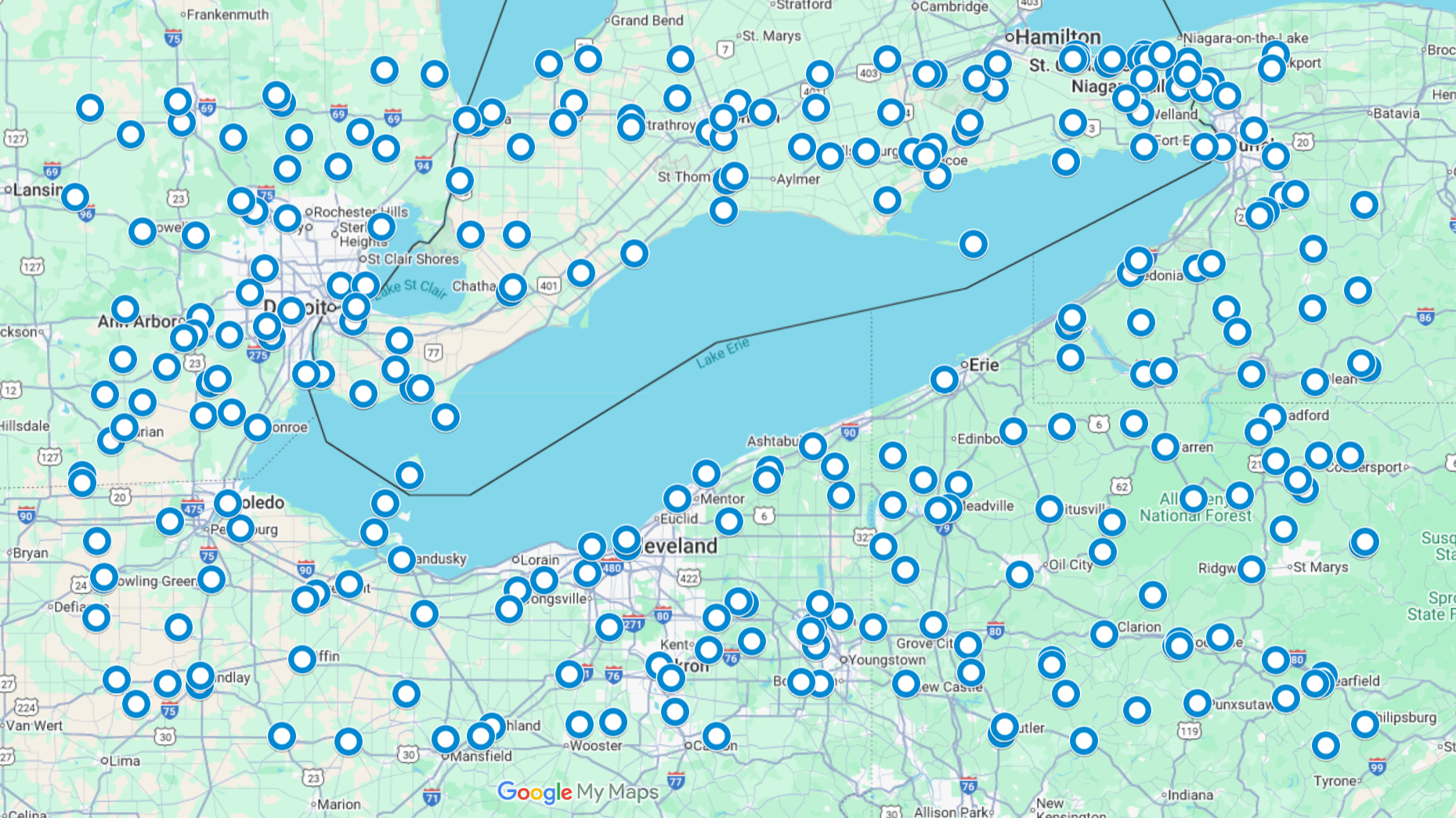}
    \end{subfigure}
    \hfill
    \begin{subfigure}[b]{0.48\textwidth}
        \centering
        \includegraphics[width=\textwidth]{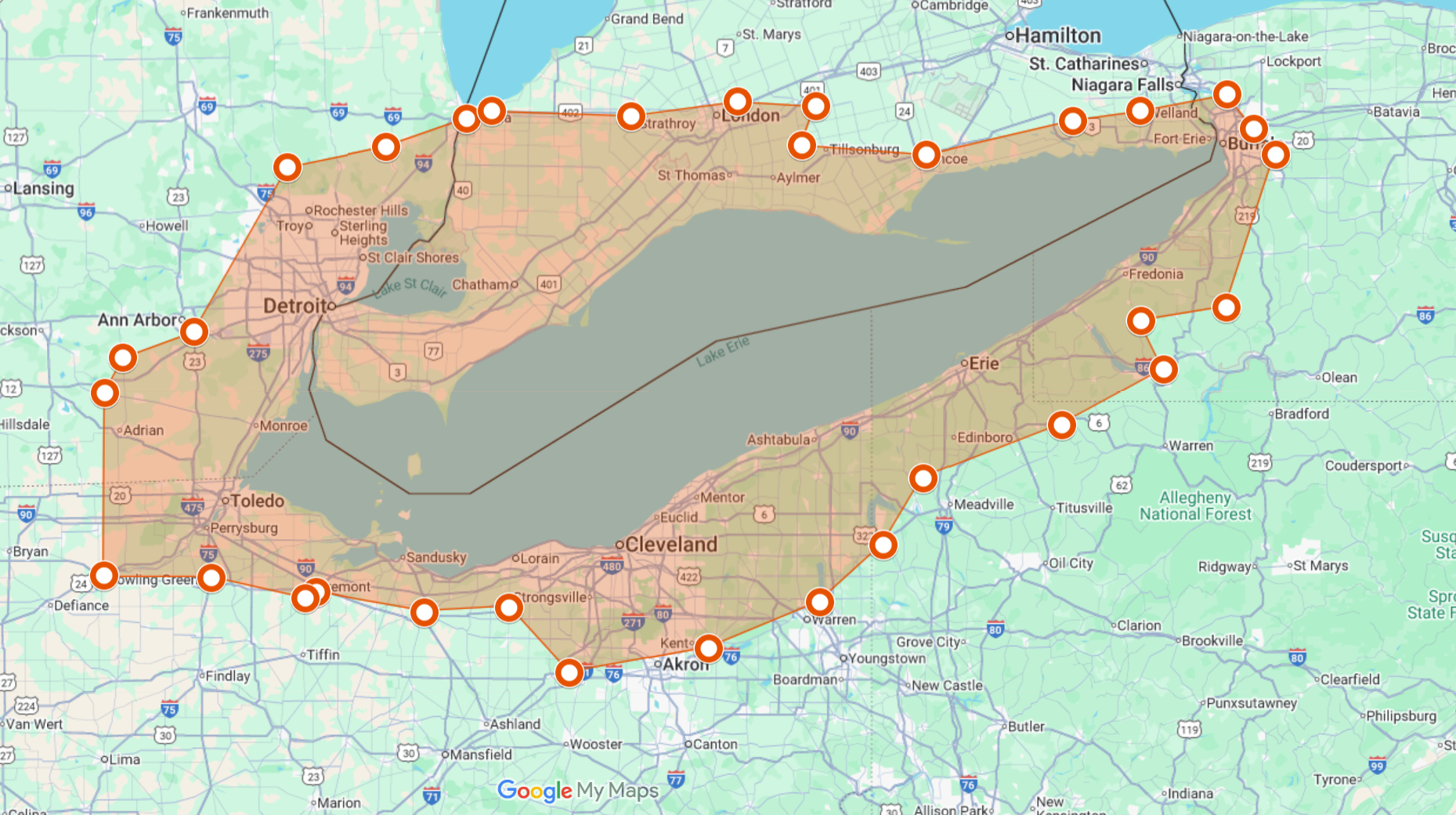}
    \end{subfigure}
    \caption{Lake Erie weather stations (left) and longest generator (right)}
    \label{fig:station_and_cyc_ON}
\end{figure}

Again, both methods detect the lake, but the barcodes tell a similar story:
\begin{figure}[H]
    \centering
    \begin{subfigure}[b]{0.48\textwidth}
        \centering
        \includegraphics[width=\textwidth]{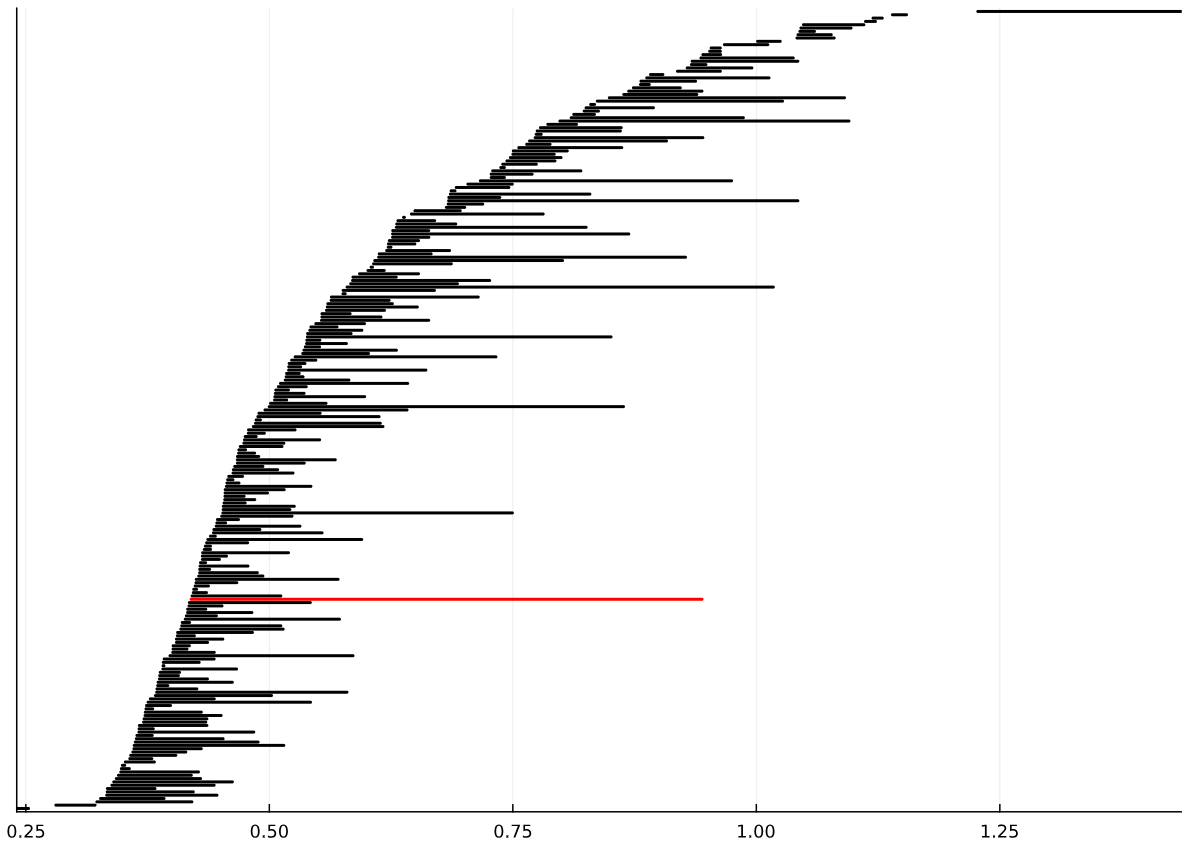}
    \end{subfigure}
    \hfill
    \begin{subfigure}[b]{0.48\textwidth}
        \centering
        \includegraphics[width=\textwidth]{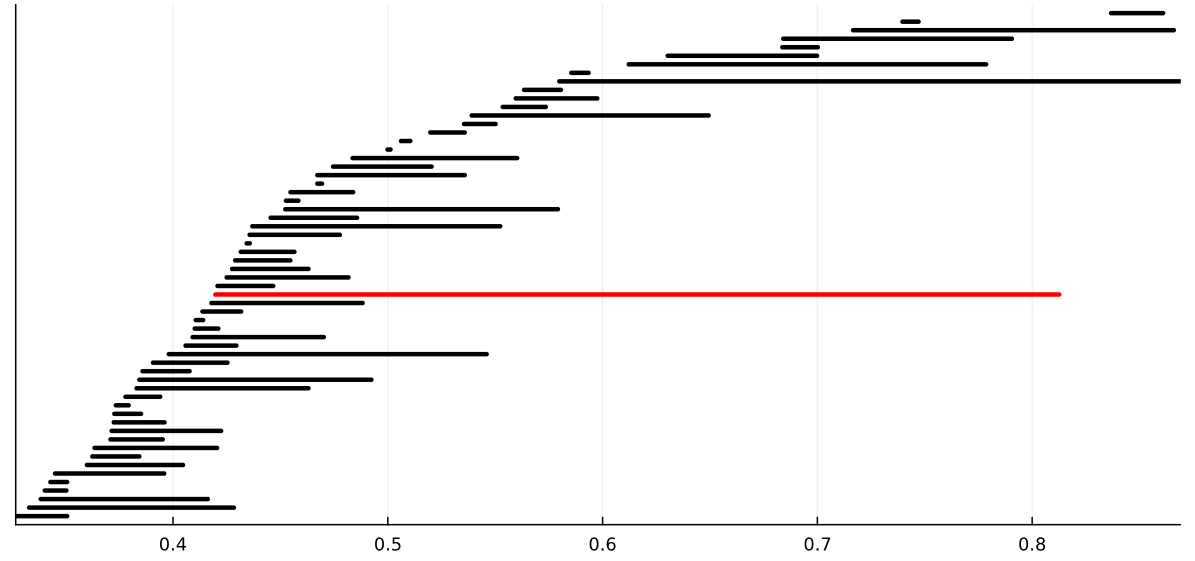}
    \end{subfigure}
    \caption{Barcodes using Vietoris-Rips (left) and discrete homology (right)}
    \label{fig:on_barcodes}
\end{figure}
\section{Summary and future directions} \label{sec:summary}
This paper introduces a new method of topological data analysis: persistent discrete homology, intended to analyze noisy data in the non-metric setting.
To ascertain persistent discrete homology's resistance to noise, we introduced the notion of overlap and used it to demonstrate that in the non-metric setting persistent discrete homology is significantly more noise resistant than Vietoris--Rips persistence.
We then compared persistent discrete homology and Vietoris--Rips persistence on realistic weather data: one given by weather simulation and one coming from NOAA.
In both cases, the outputs of persistent discrete homology were cleaner and more accurate than those of Vietoris--Rips persistence.

There are several potential directions for future work.
First, we restrict attention to degree 1 homology throughout this paper, because discrete homology is not yet efficiently computable in higher dimensions.
As more advances are made in this direction \cite{ender-kapulkin-higher}, carrying out analogous tests in higher dimensions becomes tractable.
A second direction is to establish noise resistance rigorously.
We have experimental evidence suggesting that, in the non-metric setting,
$$\mathbb{E}[O(\dgm_{DH}(X),\dgm_{DH}(X'))] > \mathbb{E}[O(\dgm_{VR}(X),\dgm_{VR}(X'))],$$
where $X'$ is obtained from $X$ by using method 2 in \cref{sec:DCH_vs_TDA};
proving this inequality for suitable setups would give a strong mathematical underpinning for the use of discrete homology.

\bibliographystyle{amsalphaurlmod}
\bibliography{all-refs}

\appendix
\renewcommand{\thesection}{\Alph{section}}

\end{document}